\documentclass[12pt,reqno]{amsart} 
\usepackage{amssymb,amsthm,amsmath}
\usepackage{hyperref}
\usepackage{float}

\theoremstyle{plain}
\newtheorem{thm}{Theorem}

\newtheorem{cor}{Corollary}

\theoremstyle{definition}
\newtheorem*{problem*}{Problem}

\newtheorem*{notation*}{Notation}

\newcommand{\be}{\begin{enumerate}}
\newcommand{\ee}{\end{enumerate}}

\title[Cyclic Subgroups of Finite Groups]{Finite groups with a prescribed number of cyclic subgroups II}
\author{Richard Belshoff}
\address{Department of Mathematics, Missouri State University, 
	Springfield, MO 65897}
\email[Corresponding author]{RBelshoff@MissouriState.edu}
\author{Joe Dillstrom}
\address{Department of Mathematics, Missouri State University, 
	Springfield, MO 65897}
\email{Dillstrom008@live.missouristate.edu}
\author{Les Reid}
\address{Department of Mathematics, Missouri State University, 
	Springfield, MO 65897}
\email{LesReid@MissouriState.edu}

\date{\today}
\begin{document}

\begin{abstract}
In \cite{tuarn}, 
T\u{a}rn\u{a}uceanu described the finite groups $G$ having exactly $|G|-1$
cyclic subgroups. In \cite{bdr},
the authors used elementary methods to completely characterize those finite groups $G$ having exactly $|G|-\Delta$ cyclic subgroups
for $\Delta=2, 3, 4$ and $5$.
In this paper, we prove that for any $\Delta >0$ if $G$ has exactly $|G|-\Delta$ cyclic subgroups, then $|G|\le 8\Delta$ and therefore the number of such $G$ is finite. We then use the computer program GAP to find all $G$ with exactly $|G|-\Delta$ cyclic subgroups for $\Delta=1,\ldots,32$.
\end{abstract}

\maketitle


\section{Introduction}
Throughout this paper, a cyclic group of order $n$ is denoted $C_n$ and 
the dihedral group of order $2n$ is denoted $D_{2n}$.
In \cite{tuarn},
Marius T\u{a}rn\u{a}uceanu proved the following theorem.
\begin{thm}[T\u{a}rn\u{a}uceanu] A finite group $G$ has exactly $|G|-1$ 
cyclic subgroups if and only if $G$ is one of the following groups:
	$C_3$, $C_4$, $S_3$, or $D_8$.
\end{thm}

Let $G$ be a finite group, and let $C(G)$ denote the poset of 
cyclic subgroups of $G$.
We define  $\Delta(G)$, or just $\Delta$ if the group is understood, to be 
the difference between the order of $G$ and the number of cyclic subgroups of $G$:
$$\Delta(G)=|G|-|C(G)|.$$  

It is known that $\Delta=0$ if and only if 
$G$ is an elementary abelian $2$-group, i.e.,
$G=C_2^n$ for some $n$.

In \cite{tuarn} where T\u{a}rn\u{a}uceanu described all groups
satisfying $\Delta=1$, he also 
posed the following natural open problem.

\begin{problem*} Describe the finite groups $G$ having $|G|-\Delta$
cyclic subgroups, where $2\le \Delta\le |G|-1$.
\end{problem*}

In \cite{bdr} we determined those finite groups for which $\Delta=2, 3, 4$, and $5$. In this paper we will show that
$|G|\le8\Delta(G)$ and use this bound and a GAP program to find all $G$ such that $\Delta(G)=1,\ldots 32$.

\section{Main Result}
Let $G$ be a finite group.  Denote by $C(G)$ the poset of cyclic subgroups of $G$,
and let $\Delta(G)=|G|-|C(G)|$. 
For any positive integer $d$, let
$n_d= |\{H\in C(G)\ :\ |H|=d\}|$ be the number of cyclic subgroups of order $d$. 
Since every element of $G$ generates a cyclic subgroup and $\phi(d)$ is
the number of generators of a cyclic group of order $d$, it follows that
$$\sum_{d \ge 1} n_d\phi(d)=|G|.$$
Because $|C(G)|=\displaystyle\sum_{d \ge 1} n_d$, we have
\begin{equation} 
\label{eqn:tuarn}
	\sum_{d\ge 1} n_d(\phi(d)-1) = \Delta(G).
	\tag{$\star$}
\end{equation}
Since $\phi(d)\ge 2$ for $d>2$, we have 

\begin{equation} 
\label{eqn:main1}
	\Delta(G) = \sum_{d\ge 1} n_d(\phi(d)-1)\ge \sum_{d>2} n_d = |C(G)| - i_2(G),
\end{equation}
where $i_2(G)$ denotes the number of elements of order 2 or less.

Now by \cite{miller},
if $G\not\cong C_2^k$, then $i_2(G)\le\frac{3}{4}|G|$ with equality if and only if $G\cong D_8\times C_2^k$ for some $k\ge 0$. If $\Delta(G)>0$, then $G\not\cong C_2^k$ so we have

\begin{eqnarray*}
	|G|&=&\Delta(G) + |C(G)|\\
 &\le& \Delta(G) + \Delta(G) + i_2(G)\\
&=&2\Delta(G)+i_2(G)\\
&\le&2\Delta(G)+\frac{3}{4}|G|.  
\end{eqnarray*}

Solving the inequality for $|G|$ yields the following result.

\begin{thm}
For any finite group $G$, if $\delta >0$ and $\Delta(G)=\delta$, then $|G|\le 8\delta$. Equality only occurs when $\delta=2^k$ and in that case $G\cong D_8\times C_2^k$.
\end{thm}

\begin{cor}
If $\delta>0$, the number of groups such that $\Delta(G)=\delta$ is finite.
\end{cor}

This gives an alternate proof of T\u{a}rn\u{a}uceanu's result. If $\Delta(G)=1$, then $|G|\le 8$ and one only need check the nine groups (up to isomorphism) of order less than 8 [$C_1,C_2,C_3,C_4,C_2\times C_2,C_5,C_6,D_6,C_7$](if $|G|=8$ we must have $G\cong D_8$).

In the remainder of this paper, we will report on the results of a GAP program that finds all groups of order less than $8\delta$ such that $\Delta(G)=\delta$ for $\delta=1,\ldots,32$.

\section{Computer Results}
Here we list all groups with difference $\delta$, for $\delta=1$ to $32$. The GAP notation 
is as follows:

\texttt{A x B} 
	denotes a direct product of a group A by a group B,

\texttt{A:B}
	denotes a split extension of A by B,

\texttt{A.B}
	denotes just an extension of A by B (split or nonsplit).

	\noindent
	For each group, we list both 
	GAP's \texttt{StructureDescription} and \texttt{GroupID}, in the
	form
$$\mbox{\texttt{StructureDescription} = \texttt{GroupID} }$$
For example, the line below that reads
$$\mbox{\texttt{D8} = \texttt[8,3] }$$
means that the dihedral group of order $8$ is the third group of order 8
in GAP's \texttt{SmallGroups} library.

\subsection*{Four groups with difference 1 (T\u{a}rn\u{a}uceanu \cite{tuarn}) }
\begin{verbatim}
C3 = [ 3, 1 ]
C4 = [ 4, 1 ]
S3 = [ 6, 1 ]
D8 = [ 8, 3 ]
\end{verbatim}

\subsection*{Four groups with difference 2 (T\u{a}rn\u{a}uceanu \cite{T2}, 
and BDR \cite{bdr})} 
\begin{verbatim}
C6 = [ 6, 2 ]
C4 x C2 = [ 8, 2 ]
D12 = [ 12, 4 ]
C2 x D8 = [ 16, 11 ]
\end{verbatim}

\subsection*{Three groups with difference 3 (BDR \cite{bdr})}
\begin{verbatim}
C5 = [ 5, 1 ]
Q8 = [ 8, 4 ]
D10 = [ 10, 1 ]
\end{verbatim}

\subsection*{Eleven groups with difference 4 (BDR \cite{bdr})}
\begin{verbatim}
C8 = [ 8, 1 ]
C3 x C3 = [ 9, 2 ]
A4 = [ 12, 3 ]
C6 x C2 = [ 12, 5 ]
(C4 x C2) : C2 = [ 16, 3 ]
D16 = [ 16, 7 ]
C4 x C2 x C2 = [ 16, 10 ]
(C4 x C2) : C2 = [ 16, 13 ]
(C3 x C3) : C2 = [ 18, 4 ]
C2 x C2 x S3 = [ 24, 14 ]
C2 x C2 x D8 = [ 32, 46 ]
\end{verbatim}

\subsection*{Three groups with difference 5 (BDR \cite{bdr})}
\begin{verbatim}
C7 = [ 7, 1 ]
C3 : C4 = [ 12, 1 ]
D14 = [ 14, 1 ]
\end{verbatim}

\subsection*{Thirteen groups with difference 6}
\begin{verbatim}
C9 = [ 9, 1 ]
C10 = [ 10, 2 ]
C12 = [ 12, 2 ]
C4 x C4 = [ 16, 2 ]
C4 : C4 = [ 16, 4 ]
QD16 = [ 16, 8 ]
C2 x Q8 = [ 16, 12 ]
D18 = [ 18, 1 ]
D20 = [ 20, 4 ]
D24 = [ 24, 6 ]
(C2 x C2 x C2 x C2) : C2 = [ 32, 27 ]
(C4 x C4) : C2 = [ 32, 34 ]
(C2 x D8) : C2 = [ 32, 49 ]
\end{verbatim}

\subsection*{Three groups with difference 7}
\begin{verbatim}
C3 x S3 = [ 18, 3 ]
(C6 x C2) : C2 = [ 24, 8 ]
S4 = [ 24, 12 ]
\end{verbatim}

\subsection*{Fifteen groups with difference 8}
\begin{verbatim}
C8 x C2 = [ 16, 5 ]
C8 : C2 = [ 16, 6 ]
Q16 = [ 16, 9 ]
C6 x C3 = [ 18, 5 ]
C5 : C4 = [ 20, 3 ]
C2 x A4 = [ 24, 13 ]
C6 x C2 x C2 = [ 24, 15 ]
C2 x ((C4 x C2) : C2) = [ 32, 22 ]
(C4 x C2 x C2) : C2 = [ 32, 28 ]
C2 x D16 = [ 32, 39 ]
C4 x C2 x C2 x C2 = [ 32, 45 ]
C2 x ((C4 x C2) : C2) = [ 32, 48 ]
C2 x ((C3 x C3) : C2) = [ 36, 13 ]
C2 x C2 x C2 x S3 = [ 48, 51 ]
C2 x C2 x C2 x D8 = [ 64, 261 ]
\end{verbatim}

\subsection*{Three groups with difference 9}
\begin{verbatim}
C11 = [ 11, 1 ]
D22 = [ 22, 1 ]
C4 x S3 = [ 24, 5 ]
\end{verbatim}

\subsection*{Eleven groups with difference 10}
\begin{verbatim}
C14 = [ 14, 2 ]
C2 x (C3 : C4) = [ 24, 7 ]
C3 x D8 = [ 24, 10 ]
D28 = [ 28, 3 ]
((C4 x C2) : C2) : C2 = [ 32, 6 ]
C4 x D8 = [ 32, 25 ]
(C4 x C2 x C2) : C2 = [ 32, 30 ]
(C4 x C4) : C2 = [ 32, 31 ]
(C2 x D8) : C2 = [ 32, 43 ]
(C2 x Q8) : C2 = [ 32, 50 ]
S3 x S3 = [ 36, 10 ]
\end{verbatim}

\subsection*{Eight groups with difference 11}
\begin{verbatim}
C13 = [ 13, 1 ]
C15 = [ 15, 1 ]
C16 = [ 16, 1 ]
C5 : C4 = [ 20, 1 ]
SL(2,3) = [ 24, 3 ]
D26 = [ 26, 1 ]
D30 = [ 30, 3 ]
D32 = [ 32, 18 ]
\end{verbatim}

\subsection*{Twenty-one groups with difference 12}
\begin{verbatim}
C18 = [ 18, 2 ]
C10 x C2 = [ 20, 5 ]
C7 : C3 = [ 21, 1 ]
C3 : Q8 = [ 24, 4 ]
C12 x C2 = [ 24, 9 ]
(C4 x C2) : C4 = [ 32, 2 ]
(C8 x C2) : C2 = [ 32, 9 ]
C4 x C4 x C2 = [ 32, 21 ]
C2 x (C4 : C4) = [ 32, 23 ]
(C4 x C4) : C2 = [ 32, 24 ]
(C2 x Q8) : C2 = [ 32, 29 ]
(C4 x C4) : C2 = [ 32, 33 ]
C2 x QD16 = [ 32, 40 ]
(C8 x C2) : C2 = [ 32, 42 ]
C2 x C2 x Q8 = [ 32, 47 ]
D36 = [ 36, 4 ]
C2 x C2 x D10 = [ 40, 13 ]
C2 x D24 = [ 48, 36 ]
C2 x ((C2 x C2 x C2 x C2) : C2) = [ 64, 202 ]
C2 x ((C4 x C4) : C2) = [ 64, 211 ]
C2 x ((C2 x D8) : C2) = [ 64, 264 ]
\end{verbatim}

\subsection*{Five groups with difference 13}
\begin{verbatim}
(C3 x C3) : C3 = [ 27, 3 ]
C3 x C3 x C3 = [ 27, 5 ]
(C3 x C3) : C4 = [ 36, 9 ]
D8 x S3 = [ 48, 38 ]
(C3 x C3 x C3) : C2 = [ 54, 14 ]
\end{verbatim}

\subsection*{Thirteen groups with difference 14}
\begin{verbatim}
C20 = [ 20, 2 ]
C3 x Q8 = [ 24, 11 ]
(C8 : C2) : C2 = [ 32, 7 ]
(C4 x C4) : C2 = [ 32, 11 ]
C4 x Q8 = [ 32, 26 ]
(C2 x C2) . (C2 x C2 x C2) = [ 32, 32 ]
C4 : Q8 = [ 32, 35 ]
(C2 x Q8) : C2 = [ 32, 44 ]
C6 x S3 = [ 36, 12 ]
D40 = [ 40, 6 ]
C2 x ((C6 x C2) : C2) = [ 48, 43 ]
C2 x S4 = [ 48, 48 ]
D8 x D8 = [ 64, 226 ]
\end{verbatim}

\subsection*{Four groups with difference 15}
\begin{verbatim}
C17 = [ 17, 1 ]
C3 : C8 = [ 24, 1 ]
QD32 = [ 32, 19 ]
D34 = [ 34, 1 ]
\end{verbatim}

\subsection*{Thirty groups with difference 16}
\begin{verbatim}
C24 = [ 24, 2 ]
C3 x D10 = [ 30, 2 ]
(C8 x C2) : C2 = [ 32, 5 ]
Q8 : C4 = [ 32, 10 ]
C8 : C4 = [ 32, 13 ]
C8 : C4 = [ 32, 14 ]
C8 x C2 x C2 = [ 32, 36 ]
C2 x (C8 : C2) = [ 32, 37 ]
(C8 x C2) : C2 = [ 32, 38 ]
C2 x Q16 = [ 32, 41 ]
C3 x A4 = [ 36, 11 ]
C6 x C6 = [ 36, 14 ]
C2 x (C5 : C4) = [ 40, 12 ]
D48 = [ 48, 7 ]
C2 x C2 x A4 = [ 48, 49 ]
(C2 x C2 x C2 x C2) : C3 = [ 48, 50 ]
C6 x C2 x C2 x C2 = [ 48, 52 ]
(C2 x ((C4 x C2) : C2)) : C2 = [ 64, 60 ]
(C2 x C2 x D8) : C2 = [ 64, 73 ]
C2 x C2 x ((C4 x C2) : C2) = [ 64, 193 ]
C2 x ((C4 x C2 x C2) : C2) = [ 64, 203 ]
(C2 x C2 x D8) : C2 = [ 64, 215 ]
(C2 x C2 x D8) : C2 = [ 64, 216 ]
C2 x C2 x D16 = [ 64, 250 ]
C4 x C2 x C2 x C2 x C2 = [ 64, 260 ]
C2 x C2 x ((C4 x C2) : C2) = [ 64, 263 ]
(C2 x ((C4 x C2) : C2)) : C2 = [ 64, 266 ]
C2 x C2 x ((C3 x C3) : C2) = [ 72, 49 ]
C2 x C2 x C2 x C2 x S3 = [ 96, 230 ]
C2 x C2 x C2 x C2 x D8 = [ 128, 2320 ]
\end{verbatim}

\subsection*{Eight groups with difference 17}
\begin{verbatim}
C19 = [ 19, 1 ]
C21 = [ 21, 2 ]
C7 : C4 = [ 28, 1 ]
(C3 x C3) : C4 = [ 36, 7 ]
D38 = [ 38, 1 ]
(C10 x C2) : C2 = [ 40, 8 ]
D42 = [ 42, 5 ]
(C4 x S3) : C2 = [ 48, 41 ]
\end{verbatim}

\subsection*{Seventeen groups with difference 18}
\begin{verbatim}
C22 = [ 22, 2 ]
C5 x C5 = [ 25, 2 ]
C8 x C4 = [ 32, 3 ]
C8 : C4 = [ 32, 4 ]
C2 . ((C4 x C2) : C2) = (C2 x C2) . (C4 x C2) = [ 32, 8 ]
C4 : C8 = [ 32, 12 ]
D44 = [ 44, 3 ]
(C12 x C2) : C2 = [ 48, 14 ]
C2 x C4 x S3 = [ 48, 35 ]
(C12 x C2) : C2 = [ 48, 37 ]
(C5 x C5) : C2 = [ 50, 4 ]
(((C4 x C2) : C2) : C2) : C2 = [ 64, 138 ]
(C8 x C4) : C2 = [ 64, 174 ]
(C2 x C2 x D8) : C2 = [ 64, 227 ]
(C4 x D8) : C2 = [ 64, 231 ]
((C4 x C2 x C2) : C2) : C2 = [ 64, 241 ]
((C4 x C4) : C2) : C2 = [ 64, 242 ]
\end{verbatim}

\subsection*{Eight groups with difference 19}
\begin{verbatim}
C9 x C3 = [ 27, 2 ]
C9 : C3 = [ 27, 4 ]
Q32 = [ 32, 20 ]
C4 x D10 = [ 40, 5 ]
(C7 : C3) : C2 = [ 42, 1 ]
(C3 x D8) : C2 = [ 48, 15 ]
(C2 x (C3 : C4)) : C2 = [ 48, 39 ]
(C9 x C3) : C2 = [ 54, 7 ]
\end{verbatim}

\subsection*{Twenty-eight groups with difference 20}
\begin{verbatim}
C14 x C2 = [ 28, 4 ]
C5 x S3 = [ 30, 1 ]
C4 . D8 = C4 . (C4 x C2) = [ 32, 15 ]
C3 x (C3 : C4) = [ 36, 6 ]
(C2 x (C3 : C4)) : C2 = [ 48, 19 ]
GL(2,3) = [ 48, 29 ]
A4 : C4 = [ 48, 30 ]
C2 x C2 x (C3 : C4) = [ 48, 42 ]
C6 x D8 = [ 48, 45 ]
C2 x C2 x D14 = [ 56, 12 ]
(C4 x C2 x C2 x C2) : C2 = [ 64, 67 ]
(C4 x C4 x C2) : C2 = [ 64, 71 ]
(C2 x ((C4 x C2) : C2)) : C2 = [ 64, 75 ]
C2 x (((C4 x C2) : C2) : C2) = [ 64, 90 ]
(C2 x C2 x D8) : C2 = [ 64, 128 ]
C2 x C4 x D8 = [ 64, 196 ]
(C4 x D8) : C2 = [ 64, 199 ]
C2 x ((C4 x C2 x C2) : C2) = [ 64, 205 ]
(C4 x C2 x C2 x C2) : C2 = [ 64, 206 ]
C2 x ((C4 x C4) : C2) = [ 64, 207 ]
(C4 x C4 x C2) : C2 = [ 64, 213 ]
(C2 x ((C4 x C2) : C2)) : C2 = [ 64, 218 ]
(C4 x D8) : C2 = [ 64, 219 ]
(C4 x D8) : C2 = [ 64, 221 ]
C2 x ((C2 x D8) : C2) = [ 64, 254 ]
(C2 x D16) : C2 = [ 64, 257 ]
C2 x ((C2 x Q8) : C2) = [ 64, 265 ]
C2 x S3 x S3 = [ 72, 46 ]
\end{verbatim}

\subsection*{Five groups with difference 21}
\begin{verbatim}
C23 = [ 23, 1 ]
C9 : C4 = [ 36, 1 ]
C12 x C3 = [ 36, 8 ]
D46 = [ 46, 1 ]
(C12 x C3) : C2 = [ 72, 33 ]
\end{verbatim}

\subsection*{Thirty groups with difference 22}
\begin{verbatim}
C25 = [ 25, 1 ]
C26 = [ 26, 2 ]
C28 = [ 28, 2 ]
C30 = [ 30, 4 ]
C16 x C2 = [ 32, 16 ]
C16 : C2 = [ 32, 17 ]
C2 x (C5 : C4) = [ 40, 7 ]
(C4 x C4) : C3 = [ 48, 3 ]
C24 : C2 = [ 48, 6 ]
C2 x SL(2,3) = [ 48, 32 ]
D50 = [ 50, 1 ]
D52 = [ 52, 4 ]
((C3 x C3) : C3) : C2 = [ 54, 5 ]
((C3 x C3) : C3) : C2 = [ 54, 8 ]
C3 x ((C3 x C3) : C2) = [ 54, 13 ]
D56 = [ 56, 5 ]
D60 = [ 60, 12 ]
(((C4 x C2) : C2) : C2) : C2 = [ 64, 34 ]
((C4 x C4) : C2) : C2 = [ 64, 134 ]
(((C4 x C2) : C2) : C2) : C2 = [ 64, 139 ]
(C4 x D8) : C2 = [ 64, 140 ]
(C2 x D16) : C2 = [ 64, 177 ]
C2 x D32 = [ 64, 186 ]
(C4 x D8) : C2 = [ 64, 228 ]
(C2 x C2 x Q8) : C2 = [ 64, 229 ]
(C4 x D8) : C2 = [ 64, 232 ]
(C4 x D8) : C2 = [ 64, 234 ]
(C4 x D8) : C2 = [ 64, 236 ]
(C4 x D8) : C2 = [ 64, 240 ]
((C2 x C2) . (C2 x C2 x C2)) : C2 = [ 64, 243 ]
\end{verbatim}

\subsection*{Four groups with difference 23}
\begin{verbatim}
C27 = [ 27, 1 ]
(C3 x Q8) : C2 = [ 48, 17 ]
Q8 x S3 = [ 48, 40 ]
D54 = [ 54, 1 ]
\end{verbatim}

\subsection*{Fifty-nine  groups with difference 24}
\begin{verbatim}
(C2 x C2) : C9 = [ 36, 3 ]
C18 x C2 = [ 36, 5 ]
C13 : C3 = [ 39, 1 ]
C5 : Q8 = [ 40, 4 ]
C10 x C2 x C2 = [ 40, 14 ]
C2 x (C7 : C3) = [ 42, 2 ]
C3 x D14 = [ 42, 4 ]
C4 x (C3 : C4) = [ 48, 11 ]
(C3 : C4) : C4 = [ 48, 12 ]
C12 : C4 = [ 48, 13 ]
C3 x ((C4 x C2) : C2) = [ 48, 21 ]
C3 x D16 = [ 48, 25 ]
C4 x A4 = [ 48, 31 ]
SL(2,3) : C2 = [ 48, 33 ]
C2 x (C3 : Q8) = [ 48, 34 ]
C12 x C2 x C2 = [ 48, 44 ]
C3 x ((C4 x C2) : C2) = [ 48, 47 ]
C13 : C4 = [ 52, 3 ]
(C4 x C2 x C2) : C4 = [ 64, 23 ]
C2 x ((C4 x C2) : C4) = [ 64, 56 ]
C4 x ((C4 x C2) : C2) = [ 64, 58 ]
(C2 x (C4 : C4)) : C2 = [ 64, 61 ]
((C4 x C2) : C4) : C2 = [ 64, 62 ]
(C2 x (C4 : C4)) : C2 = [ 64, 66 ]
(C4 x C4 x C2) : C2 = [ 64, 69 ]
(C2 x C2 x Q8) : C2 = [ 64, 74 ]
(C2 x (C4 : C4)) : C2 = [ 64, 77 ]
(C2 x (C4 : C4)) : C2 = [ 64, 78 ]
(C2 x (C4 : C4)) : C2 = [ 64, 80 ]
(((C4 x C2) : C2) : C2) : C2 = [ 64, 91 ]
C2 x ((C8 x C2) : C2) = [ 64, 95 ]
(C2 x (C8 : C2)) : C2 = [ 64, 99 ]
(C2 x D16) : C2 = [ 64, 130 ]
(C2 x QD16) : C2 = [ 64, 131 ]
(C8 x C2 x C2) : C2 = [ 64, 147 ]
(C2 x (C8 : C2)) : C2 = [ 64, 150 ]
C4 x C4 x C2 x C2 = [ 64, 192 ]
C2 x C2 x (C4 : C4) = [ 64, 194 ]
C2 x ((C4 x C4) : C2) = [ 64, 195 ]
C4 x ((C4 x C2) : C2) = [ 64, 198 ]
(C4 x Q8) : C2 = [ 64, 201 ]
C2 x ((C2 x Q8) : C2) = [ 64, 204 ]
C2 x ((C4 x C4) : C2) = [ 64, 209 ]
(C4 x C4 x C2) : C2 = [ 64, 210 ]
(C2 x C2 x Q8) : C2 = [ 64, 217 ]
(C4 x D8) : C2 = [ 64, 220 ]
(C4 x Q8) : C2 = [ 64, 223 ]
((C2 x Q8) : C2) : C2 = [ 64, 224 ]
C2 x C2 x QD16 = [ 64, 251 ]
C2 x ((C8 x C2) : C2) = [ 64, 253 ]
(C2 x (C8 : C2)) : C2 = [ 64, 256 ]
(C2 x QD16) : C2 = [ 64, 258 ]
C2 x C2 x C2 x Q8 = [ 64, 262 ]
C2 x C2 x D18 = [ 72, 17 ]
C2 x C2 x C2 x D10 = [ 80, 51 ]
C2 x C2 x D24 = [ 96, 207 ]
C2 x C2 x ((C2 x C2 x C2 x C2) : C2) = [ 128, 2163 ]
C2 x C2 x ((C4 x C4) : C2) = [ 128, 2172 ]
C2 x C2 x ((C2 x D8) : C2) = [ 128, 2323 ]
\end{verbatim}

\subsection*{Six groups with difference 25}
\begin{verbatim}
(C3 : C8) : C2 = [ 48, 16 ]
C3 x C3 x S3 = [ 54, 12 ]
S3 x D10 = [ 60, 8 ]
(C6 x C6) : C2 = [ 72, 35 ]
(S3 x S3) : C2 = [ 72, 40 ]
(C3 x A4) : C2 = [ 72, 43 ]
\end{verbatim}

\subsection*{Twenty-nine groups with difference 26}
\begin{verbatim}
C32 = [ 32, 1 ]
C5 : C8 = [ 40, 3 ]
C5 x D8 = [ 40, 10 ]
C2 . S4 = SL(2,3) . C2 = [ 48, 28 ]
C2 x ((C3 x C3) : C3) = [ 54, 10 ]
C6 x C3 x C3 = [ 54, 15 ]
((C8 : C2) : C2) : C2 = [ 64, 32 ]
(C4 x C4) : C4 = [ 64, 35 ]
D64 = [ 64, 52 ]
C4 x D16 = [ 64, 118 ]
(C4 x D8) : C2 = [ 64, 123 ]
((C4 x C4) : C2) : C2 = [ 64, 135 ]
((C4 x C4) : C2) : C2 = [ 64, 136 ]
(C2 x QD16) : C2 = [ 64, 141 ]
(C4 x D8) : C2 = [ 64, 144 ]
(C8 x C4) : C2 = [ 64, 167 ]
((C8 x C2) : C2) : C2 = [ 64, 171 ]
(C8 x C4) : C2 = [ 64, 173 ]
(C8 x C4) : C2 = [ 64, 176 ]
(C2 x D16) : C2 = [ 64, 190 ]
Q8 x D8 = [ 64, 230 ]
(C4 x Q8) : C2 = [ 64, 233 ]
(C4 x Q8) : C2 = [ 64, 235 ]
(C4 x Q8) : C2 = [ 64, 237 ]
((C4 x C4) : C2) : C2 = [ 64, 244 ]
(C6 x S3) : C2 = [ 72, 23 ]
C2 x ((C3 x C3) : C4) = [ 72, 45 ]
C2 x D8 x S3 = [ 96, 209 ]
C2 x ((C3 x C3 x C3) : C2) = [ 108, 44 ]
\end{verbatim}

\subsection*{Five groups with difference 27}
\begin{verbatim}
C29 = [ 29, 1 ]
C36 = [ 36, 2 ]
(C14 x C2) : C2 = [ 56, 7 ]
D58 = [ 58, 1 ]
D72 = [ 72, 6 ]
\end{verbatim}

\subsection*{Fifty-seven groups with difference 28}
\begin{verbatim}
C20 x C2 = [ 40, 9 ]
C8 x S3 = [ 48, 4 ]
C24 : C2 = [ 48, 5 ]
C3 : Q16 = [ 48, 8 ]
C12 x C4 = [ 48, 20 ]
C3 x (C4 : C4) = [ 48, 22 ]
C3 x QD16 = [ 48, 26 ]
C6 x Q8 = [ 48, 46 ]
C3 x D18 = [ 54, 3 ]
(C9 : C3) : C2 = [ 54, 6 ]
A5 = [ 60, 5 ]
((C8 x C2) : C2) : C2 = [ 64, 8 ]
C4 x C4 x C4 = [ 64, 55 ]
(C4 x C4) : C4 = [ 64, 57 ]
C4 x (C4 : C4) = [ 64, 59 ]
(C4 x C4) : C4 = [ 64, 63 ]
(C4 x C4) : C4 = [ 64, 64 ]
(C4 x C4) : C4 = [ 64, 65 ]
(C4 : C4) : C4 = [ 64, 68 ]
(C4 : C4) : C4 = [ 64, 70 ]
(C2 x Q8) : C4 = [ 64, 72 ]
(C4 x C2) : Q8 = [ 64, 76 ]
(C2 x C2 x C2) . (C2 x C2 x C2) = [ 64, 79 ]
(C2 x C2 x C2) . (C2 x C2 x C2) = [ 64, 81 ]
(C2 x C2 x C2) . (C2 x C2 x C2) = [ 64, 82 ]
C2 x ((C8 : C2) : C2) = [ 64, 92 ]
(C8 x C2 x C2) : C2 = [ 64, 97 ]
(C2 x (C8 : C2)) : C2 = [ 64, 98 ]
C2 x ((C4 x C4) : C2) = [ 64, 101 ]
(C2 x (C8 : C2)) : C2 = [ 64, 102 ]
(C2 x C2 x Q8) : C2 = [ 64, 129 ]
(C2 x Q16) : C2 = [ 64, 133 ]
(C8 x C2 x C2) : C2 = [ 64, 146 ]
(C2 x (C8 : C2)) : C2 = [ 64, 149 ]
(C2 x D16) : C2 = [ 64, 153 ]
(C2 x (C4 : C4)) : C2 = [ 64, 161 ]
(C2 x (C4 : C4)) : C2 = [ 64, 162 ]
((C8 x C2) : C2) : C2 = [ 64, 163 ]
C2 x C4 x Q8 = [ 64, 197 ]
(C4 x Q8) : C2 = [ 64, 200 ]
C2 x ((C2 x C2) . (C2 x C2 x C2)) = [ 64, 208 ]
C2 x (C4 : Q8) = [ 64, 212 ]
(C4 x Q8) : C2 = [ 64, 214 ]
(C4 x Q8) : C2 = [ 64, 222 ]
(C4 : Q8) : C2 = [ 64, 225 ]
C2 x ((C2 x Q8) : C2) = [ 64, 255 ]
(C2 x Q16) : C2 = [ 64, 259 ]
A4 x S3 = [ 72, 44 ]
C2 x C6 x S3 = [ 72, 48 ]
C2 x D40 = [ 80, 37 ]
(C12 x C4) : C2 = [ 96, 81 ]
C2 x C2 x ((C6 x C2) : C2) = [ 96, 219 ]
C2 x C2 x S4 = [ 96, 226 ]
(C2 x C2 x C2 x C2 x C2 x C2) : C2 = [ 128, 1578 ]
(C4 x C4 x C4) : C2 = [ 128, 1599 ]
C2 x D8 x D8 = [ 128, 2194 ]
(C2 x ((C2 x D8) : C2)) : C2 = [ 128, 2326 ]
\end{verbatim}

\subsection*{Nine groups with difference 29}
\begin{verbatim}
C31 = [ 31, 1 ]
C33 = [ 33, 1 ]
C5 : C8 = [ 40, 1 ]
C11 : C4 = [ 44, 1 ]
C3 : Q16 = [ 48, 18 ]
C4 x D14 = [ 56, 4 ]
D62 = [ 62, 1 ]
D66 = [ 66, 3 ]
(C6 x S3) : C2 = [ 72, 22 ]
\end{verbatim}

\subsection*{Twenty-seven groups with difference 30}
\begin{verbatim}
C34 = [ 34, 2 ]
C5 x Q8 = [ 40, 11 ]
C2 x (C3 : C8) = [ 48, 9 ]
(C3 : C8) : C2 = [ 48, 10 ]
(C4 : C8) : C2 = [ 64, 12 ]
(C4 x C2 x C2) : C4 = [ 64, 33 ]
(C16 x C2) : C2 = [ 64, 38 ]
C4 x QD16 = [ 64, 119 ]
(C4 x Q8) : C2 = [ 64, 121 ]
((C4 x C4) : C2) : C2 = [ 64, 137 ]
(Q8 : C4) : C2 = [ 64, 142 ]
(C2 x Q16) : C2 = [ 64, 145 ]
(C8 : C4) : C2 = [ 64, 155 ]
(C8 : C4) : C2 = [ 64, 157 ]
(C8 : C4) : C2 = [ 64, 159 ]
(C8 x C4) : C2 = [ 64, 169 ]
(Q8 : C4) : C2 = [ 64, 170 ]
(C2 x Q16) : C2 = [ 64, 178 ]
C2 x QD32 = [ 64, 187 ]
(C16 x C2) : C2 = [ 64, 189 ]
Q8 : Q8 = [ 64, 238 ]
Q8 x Q8 = [ 64, 239 ]
(C2 x C2) . (C2 x C2 x C2 x C2) = [ 64, 245 ]
D68 = [ 68, 4 ]
C4 x ((C3 x C3) : C2) = [ 72, 32 ]
(C2 x C2 x C2 x S3) : C2 = [ 96, 89 ]
(D8 x S3) : C2 = [ 96, 216 ]
\end{verbatim}

\subsection*{Five groups with difference 31}
\begin{verbatim}
C35 = [ 35, 1 ]
C15 : C4 = [ 60, 7 ]
D70 = [ 70, 3 ]
(C3 x C3) : Q8 = [ 72, 41 ]
D8 x D10 = [ 80, 39 ]
\end{verbatim}

\subsection*{Sixty-four groups with difference 32}
\begin{verbatim}
C40 = [ 40, 2 ]
C7 x S3 = [ 42, 3 ]
C24 x C2 = [ 48, 23 ]
C3 x (C8 : C2) = [ 48, 24 ]
C3 x Q16 = [ 48, 27 ]
C6 x D10 = [ 60, 10 ]
((C8 x C2) : C2) : C2 = [ 64, 4 ]
(C2 x Q8) : C4 = [ 64, 9 ]
(C8 x C2) : C4 = [ 64, 18 ]
(C4 x C4) : C4 = [ 64, 20 ]
(C8 x C2) : C4 = [ 64, 21 ]
C2 x ((C8 x C2) : C2) = [ 64, 87 ]
(C2 x (C8 : C2)) : C2 = [ 64, 88 ]
(C8 x C2 x C2) : C2 = [ 64, 89 ]
(C2 x (C8 : C2)) : C2 = [ 64, 94 ]
C2 x (Q8 : C4) = [ 64, 96 ]
(Q8 : C4) : C2 = [ 64, 100 ]
C2 x (C8 : C4) = [ 64, 106 ]
C2 x (C8 : C4) = [ 64, 107 ]
(C8 : C4) : C2 = [ 64, 108 ]
(C8 : C4) : C2 = [ 64, 109 ]
(C2 x Q16) : C2 = [ 64, 132 ]
(C2 x Q16) : C2 = [ 64, 148 ]
(C2 x Q16) : C2 = [ 64, 151 ]
(C2 x QD16) : C2 = [ 64, 152 ]
(Q8 : C4) : C2 = [ 64, 164 ]
(Q8 : C4) : C2 = [ 64, 165 ]
(C8 : C4) : C2 = [ 64, 166 ]
C8 x C2 x C2 x C2 = [ 64, 246 ]
C2 x C2 x (C8 : C2) = [ 64, 247 ]
C2 x ((C8 x C2) : C2) = [ 64, 248 ]
(C2 x (C8 : C2)) : C2 = [ 64, 249 ]
C2 x C2 x Q16 = [ 64, 252 ]
C17 : C4 = [ 68, 3 ]
(C3 x (C3 : C4)) : C2 = [ 72, 21 ]
C6 x A4 = [ 72, 47 ]
C6 x C6 x C2 = [ 72, 50 ]
D80 = [ 80, 7 ]
(C2 x (C5 : C4)) : C2 = [ 80, 34 ]
C2 x C2 x (C5 : C4) = [ 80, 50 ]
C2 x D48 = [ 96, 110 ]
(C2 x C2 x C2 x S3) : C2 = [ 96, 144 ]
(C6 x D8) : C2 = [ 96, 147 ]
(C6 x D8) : C2 = [ 96, 211 ]
C2 x C2 x C2 x A4 = [ 96, 228 ]
C2 x ((C2 x C2 x C2 x C2) : C3) = [ 96, 229 ]
C6 x C2 x C2 x C2 x C2 = [ 96, 231 ]
C2 x ((C2 x ((C4 x C2) : C2)) : C2) = [ 128, 1009 ]
C2 x ((C2 x C2 x D8) : C2) = [ 128, 1116 ]
(C2 x C2 x C2 x D8) : C2 = [ 128, 1135 ]
(C2 x C2 x C2 x D8) : C2 = [ 128, 1165 ]
C2 x C2 x C2 x ((C4 x C2) : C2) = [ 128, 2151 ]
C2 x C2 x ((C4 x C2 x C2) : C2) = [ 128, 2164 ]
C2 x ((C2 x C2 x D8) : C2) = [ 128, 2177 ]
C2 x ((C2 x C2 x D8) : C2) = [ 128, 2178 ]
(C2 x C2 x C2 x D8) : C2 = [ 128, 2216 ]
(C2 x ((C4 x C4) : C2)) : C2 = [ 128, 2230 ]
C2 x C2 x C2 x D16 = [ 128, 2306 ]
C4 x C2 x C2 x C2 x C2 x C2 = [ 128, 2319 ]
C2 x C2 x C2 x ((C4 x C2) : C2) = [ 128, 2322 ]
C2 x ((C2 x ((C4 x C2) : C2)) : C2) = [ 128, 2325 ]
C2 x C2 x C2 x ((C3 x C3) : C2) = [ 144, 196 ]
C2 x C2 x C2 x C2 x C2 x S3 = [ 192, 1542 ]
C2 x C2 x C2 x C2 x C2 x D8 = [ 256, 56083 ]
\end{verbatim}

\section{Acknowledgement}
This paper is a continuation of \cite{bdr},
which was an expansion and revision of the Master's thesis of the 
second author, directed by the first and third authors.

\end{document}